\documentclass[twoside,12pt]{article}
\usepackage{amsmath}
\usepackage[usenames,dvipsnames]{color}
\usepackage{cite}
\usepackage{amsmath}
\usepackage{amsthm}
\usepackage{amssymb}

\definecolor{rosso}{rgb}{0.8,0,0}

\def\gianni #1{{\color{blue}#1}}
\def\gian #1{{\color{cyan}#1}}

\def\pier #1{{\color{rosso} #1}}

\def\gianni #1{#1}
\def\gian #1{#1}

\def\pier #1{#1}

\newcommand{\vp}{\varphi(\alpha)}
\newcommand{\dn}{\partial_{\bf n}}
\newcommand{\pt}{\partial_t}
\newcommand{\beq}{\begin{equation}}
\newcommand{\eeq}{\end{equation}}
\newcommand{\beqa}{\begin{eqnarray}}
\newcommand{\eeqa}{\end{eqnarray}}

\newcommand{\rz}{{\rm I\!R}}
\newcommand{\nz}{{\rm I\!N}}
\newcommand{\oma}{\Omega}

\newcommand{\uad}{{\cal U}_{\rm ad}}
\newcommand{\lzo}{L^2(\Omega)}
\newcommand{\lio}{L^\infty(\Omega)}

\newcommand{\lzq}{L^2(Q)}
\newcommand{\liq}{L^\infty(Q)}
\newcommand{\txinto}{\int_0^t\!\!\!\int_\Omega}
\newcommand{\texinto}{\int_0^T\!\!\!\int_\Omega}
\newcommand{\tint}{\int_0^t}

\newcommand{\rla}{\rho^\alpha}
\newcommand{\mula}{\mu^\alpha}
\newcommand{\pla}{p^\alpha}
\newcommand{\qla}{q^\alpha}
\newcommand{\cs}{{\cal S}}

\newcommand{\lzht}{L^2(0,t;H)}
\newcommand{\lzvt}{L^2(0,t;V)}

\newcommand{\dx}{\,dx}

\newcommand{\dt}{\,dt}
\newcommand{\ds}{\,ds}

\newcommand{\cb}{{\cal B}}
\renewcommand{\qed}{\hfill \colorbox{black}{\hspace{-0.01cm}}}
\renewcommand{\min}{\mathop{\rm Min}}

\def\iO{\int_\Omega}
\def\itt{\int_t^T\!\!\!\int_\Omega}

\begin{document}

\begin{center}
{\bf {\huge Distributed optimal control\\[1mm]  of a nonstandard nonlocal \\[1mm] 
phase field system with\\[4mm] double obstacle potential\footnote{This work received \pier{a partial support from the GNAMPA} (\pier{Gruppo} Nazionale per l'Analisi
Matematica, la Probabilit\`{a} e loro Applicazioni) of INDAM (Istituto Nazionale
di Alta \pier{Matematica}) and the \pier{IMATI -- C.N.R. Pavia} for PC and GG.}}}

\vspace{9mm}
{\large Pierluigi Colli$^{\!\dagger}$,
Gianni Gilardi\footnote{Dipartimento di Matematica  ``F. Casorati'',
Universit\`a di Pavia, Via Ferrata, 1, 27100 Pavia,  Italy,
(e-mail: pierluigi.colli@unipv.it, gianni.gilardi@unipv.it) },
and \\[1mm]J\"urgen Sprekels\footnote{Weierstrass Institute for 
Applied Analysis and Stochastics,
Mohrenstra\ss e 39, 10117 Berlin and Department of Mathematics, 
Humboldt-Universit\"at zu Berlin, Unter den Linden 6, 10099 Berlin, \pier{Germany}
(e-mail: juergen.sprekels@wias-berlin.de) }}\\[6mm]{\small Key words: Distributed optimal control, \pier{phase field systems, double obstacle potentials,} nonlocal operators, 
first-order necessary optimality conditions.}
\end{center}

\vspace{7mm}
\begin{abstract} \noindent
\pier{This paper is concerned with} 
a distributed optimal control problem for a nonlocal phase field  
model of Cahn--Hilliard type, which is a nonlocal version of a model 
for two-species phase segregation on an atomic lattice under the presence of diffusion. 
\pier{The \emph{local} model} has been
\pier{investigated} in a series of papers by P. Podio-Guidugli and 
the present authors\pier{; the \emph{nonlocal} model studied here} 
consists of a highly nonlinear parabolic equation
coupled to an ordinary differential inclusion of subdifferential type. The 
inclusion originates from a free energy containing the indicator function of the
interval in which the order parameter of the phase segregation attains its values. 
It also contains a nonlocal term modeling long-range interactions. Due to the strong
nonlinear couplings between the state variables (which even involve products with time derivatives), 
the analysis of the state system is difficult. In addition, the presence of the 
differential inclusion is the reason that standard arguments of optimal control theory 
cannot be applied to guarantee the existence of Lagrange multipliers. In this paper, 
we employ recent results proved  for smooth logarithmic potentials and perform a
so-called `deep quench' approximation to establish existence and first-order necessary optimality conditions for the nonsmooth case of the double obstacle potential.
\end{abstract}


\thispagestyle{empty}
\pagestyle{myheadings}
\newcommand\testopari{\sc \pier{Colli \ --- \ Gilardi \ --- \ Sprekels}}
\newcommand\testodispari{\sc \pier{Optimal control of a phase field system with double obstacle}}
\markboth{\testodispari}{\testopari}



\section{Introduction}
Let $\Omega\subset\rz^3$ denote an open and bounded domain whose smooth
boundary $\Gamma$ has the outward unit normal ${\bf n}$, let $T>0$ be a given final time, and 
let $Q:=\oma\times (0,T)$ and
$\Sigma:=\Gamma\times (0,T)$. We study in this paper distributed 
optimal control problems of the following form:

\vspace{3mm}
$({\bf P}_0)$ \,\,Minimize the cost functional
\begin{eqnarray}
\label{cost}
 J((\mu,\rho),u)&\!\!=\!\!&\frac {\beta_1} 2 \texinto |\rho-\rho_Q|^2\dx\dt
\,+\, \frac{\beta_2}{2}\texinto |\mu-\mu_Q|^2 \,dx\,dt\nonumber\\  
&&+\,\frac{\beta_3}{2}\int_0^T\!\!\!\int_\Omega |u|^2\,dx\,dt
\end{eqnarray}
subject to the state system
\begin{align}
\label{ss1}
(1+2\,g(\rho))\,\pt\mu+\mu\,g'(\rho)\,\pt\rho-\Delta\mu=u\,\quad\mbox{a.\,e. in }\,Q,\\
\label{ss2}
\pt\rho\,+\,\xi\,+\,F'(\rho)\,+\,\cb[\rho]=\mu\,g'(\rho)\,\quad\mbox{a.\,e. in }\,Q,\\
\label{ss3}
\xi\in\partial I_{[0,1]}(\rho)\,\quad\mbox{a.\,e. in }\,Q,\\
\label{ss4}
\dn  \mu=0\,\quad\mbox{a.\,e. on }\Sigma,\\
\label{ss5}
\rho(\cdot,0)=\rho_0\,,\quad \mu(\cdot,0)=\mu_0,\,\quad\mbox{a.\,e. in } \,\Omega,
\end{align}
and to the control constraints
\begin{eqnarray}
  \label{uad}
u\in\uad&\!\!:=\!\!&\left\{u\in H^1(0,T;\lzo)\,:\,0\le u\le u_{\rm max} 
\,\,\,\,\mbox{a.\,e. in }\,Q \right.\nonumber\\
&&\quad\left.\mbox{and }\,\|u\|_{H^1(0,T;\lzo)}\,\le\,R \right\}.
\end{eqnarray}

Here, $\beta_i\ge 0$, $i=1,2,3$, and $R>0$, are given constants \pier{such that
$\beta_1 + \beta_2 + \beta_3 >0 $}, and 
the threshold function $u_{\rm max}\in\liq$ is nonnegative. 
Moreover, $\rho_Q,\mu_Q\in L^2(Q)$ represent prescribed target functions 
of the tracking-type functional~$J$. 
Although  more general cost functionals could be admitted 
for large parts of the subsequent analysis, we restrict ourselves to the above situation
for the sake of a simpler exposition.

The state system (\ref{ss1})--(\ref{ss5}) constitutes a {\em nonlocal} version of a
phase field model of Cahn--Hilliard type describing  phase segregation of 
two species (atoms and vacancies, say) on a lattice, which was recently studied in~\cite{CGS3}. 
In the (simpler) original {\em local} model, which was introduced in~\cite{PG} \pier{(see also \cite{CGPS3} and \cite{CGPS5})},
the nonlocal term ${\cal B}[\rho]$ is a replaced by the diffusive term~$-\Delta\rho$.
The local model has been intensively \gianni{discussed} in the past years 
(\pier{cf.\cite{CGKPS,CGKS1,CGKS2,CGPS3,CGPS7,CGPS6,CGPS4,CGPS5,CGSRendiconti}}). 
In particular, in \cite{CGPSco} the analogue of the control
problem $({\bf P}_0)$ for the local case was investigated for $g(\rho)=\rho$; for this special
case, also the optimal boundary control problem was studied (see \cite{CGSco1}). 

The state variables of the model are the {\em order parameter} $\rho$, interpreted as a
volumetric density, and the \emph{chemical potential} $\mu$; for physical reasons,
we must have $0\le \rho\le 1$ and $\mu > 0$ almost everywhere in $Q$. The control
function $u$ on the right-hand side of (\ref{ss1}) plays the role of a 
{\em microenergy source}. We remark at this place that the requirement encoded in the definition
of~$\uad$, namely that $u$ be nonnegative, is indispensable for the forthcoming analysis;
it is needed to guarantee the nonnegativity of the chemical potential~$\mu$.
 
The nonlinearity $F$ is assumed smooth, while $I_{[0,1]}$ is the indicator
function of the interval $[0,1]$, so that the specific local free  energy
$F_{\rm loc}:=I_{[0,1]}+F$ is typically a (nonsmooth) double obstacle potential. In this connection,
the subdifferential $\,\partial I_{[0,1]}\,$ of \,$I_{[0,1]}$\, is defined by
\begin{equation*}
\xi\in\partial I_{[0,1]}(\rho) \quad \Longleftrightarrow \quad \xi \left\{
\begin{array}{ll}
\le 0 &\quad\mbox{if $\,\pier{\rho = 0}$}\\
=0 &\quad\mbox{if $\,0<\rho<1$}\\
\ge 0&\quad\mbox{if $\, \pier{\rho = 1}$} 
\end{array}
\right. .
\end{equation*}
 
The  presence of the nonlocal term $\cb[\rho]$ in \eqref{ss2} constitutes the main 
difference to the local model. 
Simple examples are given by \gianni{integral operators of the form
\beq\label{intop1}
\cb[\rho](x,t)
= \tint\iO k(t,s,x,y)\,\rho(y,s)\,ds\,dy
\eeq
}%
and \gianni{purely} spatial convolutions \gianni{like}
\beq\label{intop2}
\cb[\rho](x,t)=\iO k(|y-x|)\,\rho(y,t)\,dy,
\eeq
with sufficiently regular kernels. 

The mathematical literature on control problems for phase field systems involving equations
of viscous or nonviscous Cahn--Hilliard type is still scarce and quite recent. We refer in this connection to the works \cite{CFGS1,CFGS2,CGS1,CGS2,HW,wn99}. Control problems
for convective Cahn--Hilliard systems were studied in \cite{RS,ZL1,ZL2}, and a few
analytical contributions were made to the coupled Cahn--Hilliard/Navier--Stokes system
(cf. \cite{FGS,FRS,HW3,HW1,HW2}). The recent contribution \cite{CGRS}
deals with the optimal control of a  Cahn--Hilliard type system arising in the modeling of
solid tumor growth. \pier{On the other hand, let us quote\cite{BBCG} for the analysis of 
the optimal control problem for a  phase field system that couples an energy balance equation with an ordinary differential inclusion characterized by the presence of the graph $\partial I_{[0,1]}$ acting on the phase variable.}

The state system (\ref{ss1})--(\ref{ss5}) is highly nonlinear, with 
nonstandard couplings involving 
time derivatives in (\ref{ss1}). Moreover, \eqref{ss2}, \eqref{ss3} define a variational
inequality, and the nonlocal term $\cb[\rho]$ is another source for possible analytical
difficulties, since the
absence of the Laplacian in \eqref{ss2} may entail a low regularity of the 
order parameter $\rho$. 
A general analysis of the state system (\ref{ss1})--(\ref{ss5}) 
was given in \cite{CGS3} for the case $u=0$ (no control), and in \cite{CGS4} we have investigated
the control problem $({\bf P}_0)$ for the case of smooth (but singular) nonlinearities.

\vspace{2mm} In this paper, we aim to employ the results established in \cite{CGS4}
to treat the nondifferentiable double obstacle case when $\xi$ satisfies the 
inclusions (\ref{ss3}). Recall that it is well known that in this case all of the classical
constraint qualifications fail, so that the existence of suitable Lagrange multipliers cannot
be guaranteed using standard methods of optimal control. Instead, our approach is guided 
by a strategy employed in \cite{CFGS1} for viscous Cahn--Hilliard systems (see also \cite{CFS} for 
for the simpler case of the Allen--Cahn equation): in~\cite{CFGS1}, necessary
optimality conditions for the double obstacle case could be established by 
performing a so-called `deep quench limit' in a family of optimal control problems
with differentiable nonlinearities of a form that had been previously treated in \cite{CGS2}.  
The general idea is briefly explained as follows: we replace the inclusion (\ref{ss3})
by
\begin{equation}
\label{110}
\xi=\varphi(\alpha)\,h'(\rho), 
\end{equation}
where $h$ is defined as the logarithmic potential
\begin{equation}
\label{defh}
h(\rho)=\rho\,\ln(\rho)\,+\,(1-\rho)\,\ln(1-\rho), \quad 0<\rho<1,\quad h(0)=h(1)=0,
\end{equation}
and  $\varphi\in C^0((0,1])$ is a positive function satisfying 
\begin{equation}
\label{defphi}
\lim_{\alpha\searrow 0}\,\varphi(\alpha)
=0.
\end{equation}
We remark that we could simply choose $\,\varphi(\alpha)
=\alpha^p\,$ for some $\,p>0$; however, there might be situations (e.\,g., in the
numerical approximation) in which it is advantageous to let 
$\,\varphi\,$  have a different behavior as $\,\alpha\searrow 0$.

\vspace{2mm} 
Now observe that $h'(\rho)=\ln\left(\frac{\rho}{1-\rho}\right)$ \,and\, 
$h''(\rho)=\frac 1 {\rho(1-\rho)}>0$\, for 
$\rho\in (0,1)$. Hence, in particular, we have that
\begin{eqnarray}
\label{eq:1.13}
&&\lim_{\alpha\searrow 0}\,\varphi(\alpha)\,h'(\rho)=0 \quad\mbox{for }\, 
0<\rho<1,\nonumber\\[2mm]
&&\lim_{\alpha\searrow 0}\Bigl(\varphi(\alpha)\,\lim_{\rho\searrow 0}h'(\rho)\Bigr)\,=\,-\infty,
\quad \lim_{\alpha\searrow 0}\Bigl(\varphi(\alpha)\,\lim_{\rho\nearrow 1}
h'(\rho)\Bigr)\,=\,+\infty\,.\qquad
\end{eqnarray}
Hence, we may regard the graphs of the functions 
$\,\varphi(\alpha)\,h'\,$ as approximations to the graph of the subdifferential
$\partial I_{[0,1]}$. 

\vspace{2mm}
We then consider, for any $\alpha\in (0,1]$, the optimal control problem (later to be denoted by $({\bf P}_\alpha)$), which results if in $({\bf P}_0)$ the relation (\ref{ss3}) is replaced by 
(\ref{110}). For this type of problem, in \cite{CGS2} the existence of optimal controls 
$u^{\alpha} \in\uad$ as well as first-order necessary optimality conditions have been derived. Proving a priori estimates (uniform in $\alpha\in (0,1]$), and 
employing compactness and monotonicity arguments, we will be able to show the following existence and approximation result: whenever $\,\{u^{\alpha_n}\}\subset\uad$ is a sequence of optimal controls for 
$({\bf P}_{\alpha_n})$, where $\alpha_n\searrow 0$ as $n\to\infty$, then there exist
a subsequence of $\{\alpha_n\}$, which is again indexed by $n$, and an optimal control 
$\overline{u}\in\uad$ of
$({\bf P}_0)$ such that

\begin{equation}
\label{eq:1.14}
u^{\alpha_n}\to\overline{u} \quad\mbox{weakly-star in ${\cal X}$ as }\,
n\to\infty,
\end{equation}
where, here and in the following,
\begin{equation}
\label{defX} 
{\cal X}:=H^1(0,T;L^2(\Omega))\cap L^\infty(Q)
\end{equation}
will always denote the control space.
In other words, optimal controls for $({\bf P}_\alpha)$ are for small $\alpha>0$ likely to be `close' to 
optimal controls for $({\bf P}_0)$. It is natural to ask if the reverse holds, i.\,e., whether every optimal control for
 $({\bf P}_0)$ can be approximated by a sequence $\,\{u^{\alpha_n}\}\,$ of optimal controls
for $({\bf P}_{\alpha_n})$, for some sequence $\alpha_n\searrow 0$. 

\vspace{2mm}
Unfortunately, we are not able to prove such a `global' result that applies to all optimal controls of $({\bf P}_0)$. However,  a `local' result can be established. To this end, let 
$\overline{u}\in\uad$ be any optimal control
for $({\bf P}_0)$. We introduce the `adapted' cost functional
\begin{equation}
\label{adapcost}
\widetilde{J}((\rho,\mu),u) \,:=\,{J}((\rho,\mu),u)\,+\,\frac 1 2
\|u-\overline{u}\|^2_{L^2(Q)}
\end{equation}
and consider for every $\alpha\in (0,1]$ the {\em adapted control problem} of minimizing $\,\widetilde{J}\,$ subject to $u\in\uad$ and to the constraint that $(\rho,\mu)$ solves the approximating system (\ref{ss1}), (\ref{ss2}), (\ref{110}), 
\pier{\eqref{ss4}, \eqref{ss5}.} It will then turn out that the following is true: 

\vspace{2mm}
(i) \,There are some sequence $\,\alpha_n\searrow 0\,$ and minimizers 
$\,{\overline{u}^{\alpha_n}}\in\uad$ of the adapted control problem 
associated with $\alpha_n$, $n\in\nz$,
such that
\begin{equation}
\label{eq:1.17}
{\overline{u}^{\alpha_n}}\to\overline{u}\quad\mbox{strongly in $L^2(Q)$
as }\,n\to \infty.
\end{equation}
(ii) It is possible to pass to the limit as $\alpha\searrow 0$ in the first-order necessary
optimality conditions corresponding to the adapted control problems associated with $\alpha\in (0,1]$ in order to derive first-order necessary optimality conditions for problem $({\bf P}_0)$.

\vspace{2mm}
The paper is organized as follows: in Section~2, we give a precise statement of the problem
under investigation, and we derive some results concerning the state system 
(\ref{ss1})--(\ref{ss5}) and 
its $\alpha$-approximation which is obtained if in $({\bf P}_0)$ the relation
 (\ref{ss3}) is replaced by the \pier{relation (\ref{110}).} 
In Section~3, we then prove the existence of optimal controls and the approximation 
result formulated above in~(i). 
The final Section~4 is devoted to the derivation of the first-order necessary 
optimality conditions, where the  strategy outlined in~(ii) is employed. 

\vspace{2mm}
Throughout this paper, we will use the following notation: \pier{for a 
(real) Banach space $X\!$, we denote} by $\,\|\cdot\|_X$ its norm and the norm of 
$X\times X\times X$, by $X'$ its dual space, and by $\langle\cdot,\cdot\rangle_X$ 
the dual pairing between $X'$ and~$X$. If $X$ is an inner product space, 
then the inner product is denoted by $(\cdot,\cdot)_X$. The only exception from this convention is given
by the $L^p$ spaces, $1\le p\le\infty$, for which we use the abbreviating notation
$\|\cdot\|_p$ for the norms in $L^p(\Omega)$. Furthermore, 
we put 
$$H:=\lzo, \quad V:=H^1(\oma), \quad W:=\{w\in H^2(\oma) :\,\dn w=0\,\,\mbox{ a.\,e. on }\Gamma\}.$$
We have the dense and continuous embeddings 
$W\subset V\subset H\cong H'\subset V'\subset W'$,
where $\langle u,v\rangle_V=(u,v)_H$ and
$\langle u,w\rangle_W=(u,w)_H$ for all $u\in H$, $v\in V$, and $w\in W$.
We will
make repeated use of the Young inequalities
\begin{align}
\label{young}
&a\,b\le \mbox{$\frac 1p|a|^p+\frac 1q|b|^q$} \quad\mbox{and} \quad
a\,b\le \gamma \,a^2 + \mbox{$\frac 1{4\gamma}$} b^2, 
\end{align} 
which are valid for all \,$a,b\in \rz$,$\gamma>0$, and $p,q\in (1,+\infty)$
with $\frac 1p+\frac 1q=1$, 
as well as of the fact that for three dimensions of space and smooth domains 
the embeddings $\,V\subset L^p(\Omega)$, $1\le p\le 6$, and 
$\,H^2(\Omega)\subset C^0(\overline{\Omega})$ are continuous and 
(in the first case only for $1\le p<6$) 
compact. In particular, there are positive constants $\widetilde K_i$, $i=1,2,3$,
 which depend only on the domain $\oma$, such that
\beqa\label{embed1}
\|v\|_6&\!\!\le\!\!&\widetilde K_1\,\|v\|_V\quad\forall\,
v\in V,\\[1mm]
\label{embed2}
\|v\,w\|_H&\!\!\le\!\!&\|v\|_6\,\|w\|_3\,\le\,
\widetilde K_2\,\|v\|_V\,\|w\|_3\quad\forall\,v,w\in V, \qquad
\\[1mm]
\label{embed3}
\|v\|_\infty&\!\!\le\!\!&\widetilde K_3\,\|v\|_{H^2(\oma)}\quad\forall\,
v\in H^2(\oma).
\eeqa
For convenience, we also put, for $t\in[0,T]$,
\beq
  Q_t := \Omega \times (0,t)
  \quad\mbox{and}\quad
  Q^t := \Omega \times (t,T).
  \label{defQt}
\eeq
\pier{Please note the difference in the position of $t$. About time derivatives of a time-dependent function $v$, we warn the reader that we will use both the notations  $\pt v, \, \pt^2 v $ and the shorter ones $v_t, \, v_{tt} $.}


\section{Problem statement and results\\ for the state system}
\setcounter{equation}{0}
Consider the optimal control problem
(\ref{cost})--(\ref{uad}). 
We make the following assumptions on the data:

\vspace{5mm}
{\bf (A1)} \,\,$F\in C^3[0,1]$,  and $g\in C^3[0,1]$ satisfies $g(\rho)\ge 0$ and $g''(\rho)\le 0$ for all $\rho\in [0,1]$.

\vspace{3mm}
{\bf (A2)} \,\,$\rho_0\in V$, $\,\mu_0\in W$, where
 $\,\mu_0\ge 0$\, a.\,e. in $\oma$, and
\begin{equation}
\label{2.1}
\inf\,\{\rho_0(x): \,x\in\Omega\}>0, \quad \sup\,\{\rho_0(x):\, x\in\Omega\}<1\,. 
\end{equation}

\vspace{3mm}
{\bf (A3)} \,\,The nonlocal operator $\cb\colon L^1(Q)\to L^1(Q)$ satisfies the following conditions:
\vspace{2mm}
{\bf (i)} \,\,\,\,For every $t\in (0,T]$, we have
\begin{equation}\label{B1}
\cb[v]|_{Q_t}=\cb[w]|_{Q_t} \,\,\mbox{ whenever }\, v|_{Q_t}=w|_{Q_t}.
\end{equation}
{\bf (ii)} \,\,\,
For all $p\in [2,+\infty]$, we have $\cb(L^p(Q_t))\subset L^p(Q_t)$ and
\beq
\label{B2}
 \|\cb[v]\|_{L^p(Q_t)}\,
\le C_{\cb,p}\left(1+\|v\|_{L^p(Q_t)}\right)
\eeq
for every $v\in L^p(Q)$ and $t\in(0,T]$.

{\bf (iii)} \,\,For every $v,w\in L^1(0,T;H)$ and $t\in (0,T]$, it holds that
\beq
\label{Bsechs} 
\tint \|\cb[v](s)-\cb[w](s)\|_6\ds\,\le\,C_{\cb}\tint\|v(s)-w(s)\|_H\ds\,.
\eeq
{\bf (iv)} \,\,\,It holds, for every
$v\in L^2(0,T;V)$ and $t\in (0,T]$, that
\begin{align}
\label{B3}
&\|\nabla \cb[v]\|_{\lzht}\,\le\,C_\cb\, \pier{\bigl( 1 + \,\|v\|_{\lzvt}\bigr)}.
\end{align}

\vspace{2mm}
{\bf (v)} \,\,\,\,For every $v\in H^1(0,T; H)$, we have $\,\pt \cb[v]\in L^2(Q)$
and
\beq
  \|\pt \cb[v]\|_{L^2(Q)} \leq C_\cb \, \pier{\bigl( 1 + \|\pt v\|_{L^2(Q)}\bigr)}.
  \label{B4}
\eeq

{\bf (vi)} \,\,\,$\cb$ is continuously Fr\'echet differentiable as a mapping
from $L^2(Q)$ into $L^2(Q)$, and the Fr\'echet derivative $D\cb[\overline{v}]\in 
{\cal L}(L^2(Q),L^2(Q))$ of $\cb$ at $\overline{v}$ has for every $\overline{v}\in L^2(Q)$ 
and $t\in(0,T]$ the following properties: 
\begin{align}
\label{B5}
&{\|D\cb[\overline{v}](w)\|_{L^p(Q_t)}\,\le\,C_\cb\,\|w\|_{L^p(Q_t)} 
\quad\forall\, w\in L^p(Q),\quad \forall\, p\in[2,6],}\\[4mm]
\label{B6}
&\|\nabla(D\cb[\overline{v}](w))\|_{L^2(Q_t)} \,\le\,C_\cb\,\|w\|_{L^2(0,t;V)}\quad\forall \,
w\in L^2(0,T;V) \pier{.}
\end{align}

\vspace*{2mm}
In the above formulas, $C_{\cb,p}$ and $C_\cb$ denote given positive structural constants. 
We also notice that \eqref{B5} implicitly requires that
$D\cb[\overline v](w)|_{Q_t}$ depends only on~$w|_{Q_t}$; this is, however,
 a consequence of~\eqref{B1}.

\medskip
The statements related to the control problem $({\bf P}_0)$ depend on the assumptions made in the introduction.
We recall them here:

\smallskip
{\bf (A4)}  \,\, $J$ and $\uad$ are defined by \eqref{cost} and \eqref{uad},
respectively, where
\begin{align}
\label{constCP}
\beta_i\ge 0\,,\,\,\,\,i=1,2,3,  \quad\pier{\beta_1+ \beta_2 + \beta_3 > 0 ,  \quad  R>0,}\ \
\\[2mm]
\label{functCP}
\rho_Q,\mu_Q\in L^2(Q), \,\,\mbox{ and }\,\, u_{\rm max} \in \liq \,\,\mbox{ is such that}\nonumber\\
u_{\rm max} \geq 0  \,\,\mbox{ a.\,e. in }\,\,Q\quad\mbox{and}\quad \uad\not=\emptyset.
\end{align}

\vspace{3mm}
{\sc Remark 1:} \,\,In \cite{CGS4} it was shown that 
the integral operator \eqref{intop2} satisfies all of the conditions \eqref{B1}--\pier{\eqref{B6}}
provided that the integral
kernel satisfies $k\in C^1(0,+\infty)$ and fulfills, with
suitable constants $C_1>0$, $C_2>0$, $0<\alpha<\frac 32$, $0<\beta<\frac 52$,  the growth conditions
$$
|k(r)|\le C_1\,r^{-\alpha}, \quad |k'(r)|\le C_2\,r^{-\beta}, \quad\forall\, r>0\,.
$$
 Notice that these growth conditions are 
met by, e.\,g., the three-dimensional Newtonian potential, where $k(r)=\frac c r$ with some
$c\not =0$.  

\vspace{5mm}
The following well-posedness result for the state system is a direct consequence 
of\pier{\cite[Thms. 2.1 and 2.2]{CGS3}}.

\vspace{5mm}
{\sc Theorem 2.1.} \quad {\em Suppose that} {\bf (A1)}--{\bf (A4)} {\em are 
satisfied. Then the state system} (\ref{ss1})--(\ref{ss5}) {\em has for every $u\in\uad$ a unique solution}
$(\rho,\mu,\xi)$ {\em such that}
\begin{align}
\label{reguss1}
&\mu\in H^1(0,T;H)\cap L^\infty(0,T;V)\cap L^2(0,T;W^{2,3/2}(\oma))\cap L^\infty(Q),\\[1mm]
\label{reguss2}
&\rho\in L^\infty(0,T;V), \quad \rho_t\in L^6(Q),\quad \xi\in L^6(Q),\\[1mm]
\label{reguss3}
&\mu\ge 0, \quad 0\le\rho\le 1,\quad \mbox{a.\,e. in }\,Q.
\end{align}

\vspace{5mm}
It follows from Theorem 2.1 that under the given assumptions the control-to-state operator
for problem $({\bf P}_0)$ is well defined on $\uad$. In the following, we denote by $\cs_0$
the operator 
\beq
\label{defs0}
\cs_0:\uad \ni u \mapsto \cs_0(u):=(\mu,\rho).
\eeq

We now turn our interest to the $\alpha$--approximating state system, which results when \eqref{ss3}
is replaced by \eqref{110}. For later convenience, we write this system explicitly:
\begin{align}
\label{sa1}
(1+2\,g(\rla))\,\pt\mula+\mula\,g'(\rla)\,\pt\rla-\Delta\mula=u\,\quad\mbox{a.\,e. in }\,Q,\\
\label{sa2}
\pt\rla\,+\,\vp h'(\rla)\,+\,F'(\rla)\,+\,\cb[\rla]=\mula\,g'(\rla)\,\quad\mbox{a.\,e. in }\,Q,\\
\label{sa3}
\dn  \mula=0\,\quad\mbox{a.\,e. on }\Sigma,\\
\label{sa4}
\rla(\cdot,0)=\rho_0\,,\quad \mula(\cdot,0)=\mu_0,\,\quad\mbox{a.\,e. in } \,\Omega.
\end{align}

For this system, we have the following regularity result (cf. \cite[Thm. 2.1]{CGS4}):

\vspace{5mm} 
{\sc Theorem 2.2.} \quad {\em Suppose that} {\bf (A1)}--{\bf (A4)} {\em are 
satisfied. Then the system} (\ref{sa1})--(\ref{sa4}) {\em has for every $u\in\uad$ and 
$\alpha\in (0,1]$ a unique
solution $(\mula,\rla)$ such that}
\begin{align}
\label{regusa1}
&\mula\in W^{1,\infty}(0,T;H)\cap H^1(0,T;V)\cap L^\infty(0,T;W)\cap C^0(\overline{Q}),\qquad\\[1mm]
\label{regusa2}
&\rla\in H^2(0,T;H)\cap \pier{W^{1,\infty}(0,T;\lio)\cap H^1(0,T;V)}.
\end{align}   
{\em Moreover, there are constants $0<\rho_*(\alpha)<\rho^*(\alpha)<1$, $\mu^*(\alpha)>0$, 
which depend only on the given data and $\alpha$, such that for every $u\in\uad$ and every $\alpha
\in (0,1]$ the 
corresponding solution $(\mula,\rla)$ satisfies}
\beq
\label{sabounds1}
0<\rho_*(\alpha)\le \rla\le\rho^*(\alpha)<1,\,\quad 0\le\mula\le \mu^*(\alpha), 
\,\quad\mbox{a.\,e. in }\,Q.
\eeq

\vspace{5mm}
From Theorem 2.2, we infer that the $\alpha$--control-to-state operator
\beq
\label{defSal}
\cs_\alpha:\uad\ni u\mapsto \cs_\alpha(u):=(\mula,\rla)
\eeq
is well defined on $\uad$ for every $\alpha\in (0,1]$. 
We now establish some  a priori estimates for $(\mula,\rla)$, which are uniform in $\alpha$. At
first, we conclude from \eqref{sabounds1}, \eqref{B2}, and {\bf (A1)}, that there is some constant $K^*_1>0$,
which depends only on the data of the problem, such that
\beq
\label{sabounds2}
\max_{0\le i\le 3} \left(\left\|F^{(i)}(\rla)\right\|_{\liq}\,+\left\|g^{(i)}(\rla)\right\|
_{\liq}\right)+\left\|\cb[\rla]\right\|_{\liq}\,\le\,K_1^*,
\eeq 
whenever $\,(\mula,\rla)=\cs_\alpha(u)\,$ for some $u\in\uad$ and some $\alpha\in (0,1]$.

Unfortunately, we cannot exclude the possibility that $\,\rho_*(\alpha)\searrow 0\,$ and/or 
$\,\rho^*(\alpha)\nearrow 1\,$ as $\,\alpha\searrow 0$. This is the main reason why
the following a priori estimate, which will however be fundamental for our analysis, is
comparatively weak.

\vspace{5mm}
{\sc Theorem 2.3.} \quad
{\em Suppose that} {\bf (A1)}--{\bf (A4)} {\em are satisfied. Then there is some constant $K^*_2>0$, which depends only on the given data, such that
for every $u\in\uad$ and every $\alpha\in (0,1]$ the
corresponding solution $(\mula,\rla)$ to} \eqref{sa1}--\eqref{sa4} {\em satisfies}
\begin{eqnarray}
\label{sabounds3}
&&\|\mula\|_
{H^1(0,T;H)\cap L^\infty(0,T;V)\cap L^2(0,T;H^2(\oma))\cap
  L^\infty(Q)}\nonumber\\[2mm]
&&+\,\|\rla\|_{
W^{1,6}(0,T;L^6(\pier{\Omega}))\pier{\cap L^\infty(0,T;V)}}\,+\,\left\|\vp\,h'(\rla)\right\|_{L^6(Q)}\,\nonumber\\[1mm]
&&+
\texinto \vp \mbox{$\,\frac 1{\rla (1-\rla)}$}  
\left|\nabla\rla\right|^2   
\dx\dt
\,\le\,K^*_2\,.
\end{eqnarray}

\vspace{3mm}
{\sc Proof:} \,\,\,Let $u\in\uad$ be fixed and $(\mula,\rla)=\cs_\alpha(u)$ for $\alpha\in (0,1]$.
We establish the validity of \eqref{sabounds3} in a series of steps. In what follows,
we denote by $C_i$, $i\in\nz$, positive constants that may depend on the data of our problem but not on
$\alpha\in (0,1]$. For the sake of a better readability, we will omit the arguments of functions if there
is no danger of confusion, as well as the superscript $\alpha\in (0,1]$, which will only be written at the
end of each step of estimation.

\vspace{2mm}
\underline{\sc Step 1:} \,\,\,First observe that 
$\pt\left((\mbox{$\frac 12$}+g(\rho))\mu^2\right)=(1+2g(\rho))\,\mu\,\mu_t + g'(\rho)\,\rho_t\,\mu^2\,.$
Hence, multiplication of \eqref{sa1} by $\,\mu\,$ and integration over $\,Q_t$, for $0<t\le T$, yields the
identity
\begin{align*}
& \iO\left(\mbox{$\frac 12$}+g(\rho(t)\right)\mu^2(t)\dx \,+\,\txinto |\nabla\mu|^2\dx\ds\\[1mm]
&=\iO\left(\mbox{$\frac 12$}+g(\rho_0)\right)\mu_0^2\dx\,+\,\txinto u\,\mu\dx\ds\,,
\end{align*}
whence, using {\bf (A1)}, Young's inequality, and Gronwall's lemma, we easily deduce that
\beq
\label{est1}
\left\|\mula\right\|_{L^\infty(0,T;H)\cap L^2(0,T;V)}\,\le\,C_1\quad\forall\,\alpha\in (0,1].
\eeq
Moreover, since the embedding $\,(L^\infty(0,T;H)\cap L^2(0,T;V))\subset L^{10/3}(Q)\,$ is continuous, we
have (by possibly choosing a larger constant $C_1$) that
\beq
\label{est11}
\left\|\mula\right\|_{L^{10/3}(Q)}\,\le\,C_1\quad\forall\,\alpha\in (0,1].
\eeq

\vspace{2mm}
\underline{\sc Step 2:} \,\,\,Next, we multiply \eqref{sa2} by $\,\rho_t\,$ and integrate over $Q_t$,
where $0<t\le T$. It follows
\begin{align*}
&\txinto |\rho_t|^2\dx\ds\,+\iO \vp \,h(\rho(t))\dx\,-\,\iO\vp\,h(\rho_0)\dx\\[1mm]
&=\,\txinto\left(-F'(\rho)-\cb[\rho]\,+\,\mu\,g'(\rho)\right)\rho_t\dx\ds\,.
\end{align*}
Thus, employing \eqref{2.1}, \eqref{sabounds2}, and \eqref{est1}, we readily obtain from Young's inequality that
\beq
\label{est2}
\left\|\rla\right\|_{H^1(0,T;H)}\,\le\,C_2 \quad\forall\,\alpha\in (0,1].
\eeq

\vspace{2mm}
\underline{\sc Step 3:} \,\,\,Now observe that we have the identity
$$
\nabla\rho_t\,+\,\vp\,h''(\rho)\nabla\rho\,=\,-F''(\rho)\nabla\rho-\nabla\cb[\rho]+g'(\rho)\nabla\mu
+\mu\,g''(\rho)\nabla\rho\,.
$$ 
Taking the euclidean scalar product of this identity with $\,\nabla\rho$, and integrating the result over
$\,Q_t$, where $0<t\le T$, we find that
\begin{align*}
&\frac 12\,\|\nabla\rho(t)\|_H^2\,-\,\frac 12\,\|\nabla\rho_0\|_H^2\,+\,\txinto\vp\,h''(\rho)\,|\nabla\rho|^2\dx\ds
\\[1mm]
&=\,\txinto\left(-F''(\rho)\nabla\rho-\nabla\cb[\rho]+g'(\rho)\nabla\mu\right)\cdot\nabla\rho\dx\ds\\[1mm]
&\quad +\,\txinto\mu\,g''(\rho)\,|\nabla\rho|^2\dx\ds\,,
\end{align*}
where, since $h''(\rho)\ge 0$, the integral on the left-hand side is nonnegative, while,
thanks to the fact that $g''(\rho)\le 0$ and $\mu\ge 0$, the second integral on the right-hand side is nonpositive.
Moreover, by using \eqref{B3}, \eqref{sabounds2}, \eqref{est1}, and Young's
inequality, we easily find that the first integral on the right-hand side is bounded by an expression of the
form $\,\,C_3\,(1+\|\rho\|_{\lzvt}^2)$. Hence, we can infer from Gronwall's lemma and \eqref{est2} that
\beq
\label{est3}
\left\|\rla\right\|_{L^\infty(0,T;V)}\,+\texinto\vp\,h''(\rla)\left|\nabla\rla\right|^2\dx\ds\,\le\,C_4
\quad\forall\,\alpha\in (0,1].
\eeq
\gianni{In particular, this yields the desired bound for the last term on the right-hand side of~\eqref{sabounds3}.}

\vspace{2mm}
\underline{\sc Step 4:} \,\,\,Next, we multiply \eqref{sa2} by $|\vp\,h'(\rho)|^{7/3}\,{\rm sign}(\vp\,h'(\rho))$
(where we put ${\rm sign}(0):=0$) and integrate over $Q_t$, where $0<t\le T$. We obtain the inequality
\begin{align}
\label{est4}
&\iO (\Phi(x,t)-\Phi(x,0))\dx \,+\,\|\vp\,h'(\rho)\|_{L^{10/3}(Q_t)}^{10/3}\nonumber\\[1mm]
&\le\,\txinto |{-}F'(\rho)-\cb[\rho]+\mu\,g'(\rho)|\,|\vp\,h'(\rho)|^{7/3}\dx\ds\,,
\end{align}
where, by the definition of~$h$, the function
$$\Phi(x,t):=\int_{1/2}^{\rho(x,t)} |\vp\,h'(r)|^{7/3}\,{\rm sign}(\vp\,h'(r))\,dr$$
is nonnegative almost everywhere in $Q$, and where $\,\Phi(\cdot,0)\,$ 
\gianni{is bounded in $\oma$ independently of $\alpha\in (0,1]$  by~\eqref{2.1}}. 
Moreover, by virtue of \eqref{sabounds2} and \eqref{est11}, we have 
$$\|{-}F'(\rho)-\cb[\rho]+\mu\,g'(\rho)\|_{L^{10/3}(Q)}\,\le\,C_5\,.$$
Employing Young's inequality with $q=10/7$, we thus can infer from \eqref{est4} that
\begin{align*}
&\|\vp\,h'(\rho)\|_{L^{10/3}(Q_t)}^{10/3}\,\le\,C_6\,+\,\frac 7{10}\,\|\vp\,h'(\rho)\|_{L^{10/3}(Q_t)}^{10/3},
\end{align*}
whence we conclude that
$$\left\|\vp\,h'(\rla)\right\|_{L^{10/3}(Q)} \,\le\,C_7\quad\forall\,\alpha\in (0,1].$$
Comparison in \eqref{sa2} then shows that also
\beq
\label{est5}
\left\|\rho_t^\alpha\right\|_{L^{10/3}(Q)}\,\le\,C_8 \quad\forall\,\alpha\in (0,1].
\eeq

\vspace{2mm}
\underline{\sc Step 5:} \,\,\,Now that \eqref{est5} is established, we are in a position to test
\eqref{sa1} by~$\mu_t$. Integrating over $Q_t$, where $0<t\le T$, and using the fact that $g$ is
nonnegative, we obtain the estimate
\begin{align}
\label{est6}
&\txinto |\mu_t|^2\dx\ds\,+\,\frac 12\,\|\nabla\mu(t)\|_H^2\,-\,\frac 12\,\|\nabla\mu_0\|_{\gianni H}^2\nonumber\\[1mm]
&\le\txinto|u|\,|\mu_t|\dx\ds\,+\,\txinto|g'(\rho)|\,|\rho_t|\,|\mu|\,|\mu_t|\dx\ds \,\,=:\,\,I_1+I_2\,.
\end{align}
Owing to Young's inequality, we have that
$$
I_1\,\le\,\frac 14\txinto|\mu_t|^2\dx\ds\,+\,C_9,
$$
and invoking \eqref{sabounds2} and \eqref{embed1}, as well as H\"older's and Young's inequalites, we find that
\begin{align*}
I_2&\le\,C_{10}\tint \|\mu_t(s)\|_2\,\|\rho_t(s)\|_3\,\|\mu(s)\|_6\ds\nonumber\\[1mm]
&\le\,\frac 14\txinto|\mu_t|^2\dx\ds\,+\,C_{11}\tint\|\rho_t(s)\|_3^2\,\|\mu(s)\|_V^2\ds\,.
\end{align*}
Now observe that $\frac {10}3>3\gianni{{}>2}$, 
and therefore H\"older's \pier{and Young's inequalities} can be employed 
to bound the $L^1(0,T)$ norm of the function 
$\,s\mapsto \|\rho_t(s)\|_3^2\,$ as follows:
\gianni{%
$$
\int_0^T\|\rho_t(s)\|_3^2\ds\,
\leq C_{12}\int_0^T \|\rho_t(s)\|_{10/3}^2\ds
\le\,C_{13}\left(1 + \int_0^T \|\rho_t(s)\|_{10/3}^{10/3}\ds\right)\,,
$$
}%
where the last expression is by \eqref{est5} bounded independently of $\alpha\in (0,1]$.
In conclusion, we may employ Gronwall's lemma to conclude from \eqref{est6} the estimate
\beq
\label{est7}
\left\|\mula\right\|_{H^1(0,T;H)\cap L^\infty(0,T;V)}\,\le\,C_{14}.
\eeq

\vspace{2mm}
\underline{\sc Step 6:} \,\,\,Next, recalling the continuity of the embedding $H^1(\oma)\subset L^6(\oma)$
and invoking \eqref{sabounds2} and \eqref{est7}, we find that
$$
\|-F'(\rho)-\cb[\rho]+\mu\,g'(\rho)\|_{L^6(Q)}\,\le C_{15}.
$$
Therefore, if we multiply \eqref{sa2} by \,\,$(\vp\,h'(\rho))^5$ and integrate over $Q_t$, where $0<t\le T$,
then we can conclude from Young's inequality with $q=6/5$ the inequality
\begin{align}
\label{est8}
&\iO(\Psi(x,t)-\Psi(x,0))\dx \,+\,\|\vp\,h'(\rho)\|_{L^6(Q_t)}^6\nonumber\\[1mm]
&\le\,C_{16}\,+\,\frac 56\,
\|\vp\,h'(\rho)\|_{L^6(Q_t)}^6\,,
\end{align}
where the function
$$
\Psi(x,t):=\int_{1/2}^{\rho(x,t)}(\vp\,h'(r))^5\,dr
$$
is obviously nonnegative,  and where $\,\Psi(\cdot,0)\,$ is by \eqref{2.1} bounded, uniformly in $\alpha\in
(0,1]$. We thus conclude from \eqref{est8} that
\beq
\label{est9}
\left\|\vp\,h'(\rla)\right\|_{L^6(Q)}\,\le\,C_{17} \quad\forall\,\alpha\in (0,1],
\eeq 
whence, by comparison in \eqref{sa2}, also
\beq
\label{est10}
\left\|\rho_t^\alpha\right\|_{L^6(Q)}\,\le\,C_{18} \quad\forall\,\alpha\in (0,1].
\eeq

\vspace{2mm}
\underline{\sc Step 7:} \,\,\,From \eqref{est10}, we have, in particular, that
\beq
\label{est12}
\left\|\rho_t^\alpha\right\|_{L^{7/3}(0,T;L^{14/3}(\oma))}\,\le\,C_{19} \quad\forall\,\alpha\in (0,1].
\eeq
We may therefore account for the conditions $\mu_0\in W$ and $\,\|u\|_{\liq}\,\le\,R\,$ in order to repeat the
argument developed in the proof of \cite[Thm.~2.3]{CGPS3}, which is based on the above summability of
$\,\rho_t^\alpha$. We should remark that the quoted proof is performed for $g(\rho)=\rho$ and $u\equiv 0$,
but only minor changes are needed to arrive at the same conclusion in the present situation (see also
the proof of the analogous \cite[Thm.~3.7]{CGPS5} in an even more complicated case). We thus can infer
the bound
\beq
\left\|\mula\right\|_{\liq}\,\le\,C_{20}\quad\forall\,\alpha\in (0,1].
\eeq

Finally, we use the previously shown estimates
to conclude that the expression 
$$\pier{(1+2\,g(\rla))\,\mu_t^\alpha +\mula\,g'(\rla)\,\rho_t^\alpha-u \ \hbox{ is bounded in }
\, \lzq }$$ 
independently of $\alpha\in (0,T]$, which then also holds true for $\,\Delta\mula$, by comparison
in \eqref{sa1}. We thus have
\begin{align}
\left\|\mula\right\|_{\gianni{L^2(0,T;H^2(\oma))}}\,\le\,C_{21}\quad\forall\,\alpha\in (0,1]. 
\end{align}
This concludes the proof of the assertion.\qed

\section{Existence and approximation of \gianni{optimal \allowbreak controls}}
\setcounter{equation}{0}
\gianni{%
Our first aim is to show the existence result stated below.
For its proof we do not use the direct method
since this would force us to find bounds for the states 
corresponding to the chosen minimizing sequence of controls. 
On the contrary, Theorem~2.3 already provides a number of estimates.
For that reason, we pass through the related approximating control problem.}

\vspace{5mm}
{\sc Theorem 3.1.} \quad {\em Suppose that the conditions} {\bf (A1)}--{\bf (A4)} {\em are fulfilled.
Then the optimal control problem} $({\bf P}_0)$ {\em admits a solution.} 

\vspace{3mm}
{\sc Proof:} \,\,\,Recall the properties \eqref{defh} and \eqref{defphi} of the functions $h$ and $\varphi$.
We consider for $\alpha\in (0,1]$ the approximating optimal control problem 

\vspace{3mm}
$({\bf P}_\alpha)$ \quad Minimize $\,J((\mula,\rla),u)\,$ subject to $u\in\uad$ and $(\mula,\rla)=\cs_\alpha(u)$,
where $\cs_\alpha$ is the solution operator of \eqref{sa1}--\eqref{sa4} defined in \eqref{defSal}.

\vspace{3mm} According to \cite[Thm.~4.1]{CGS4}, $({\bf P}_\alpha)$ enjoys for every $\alpha\in (0,1]$ at least one
solution~\gianni{$u^\alpha$}. 
Now let $\{\alpha_n\}\subset (0,1]$ be any sequence such that $\alpha_n\searrow 0$ as $n\to\infty$. By 
virtue of Theorem~2.3, \pier{there exists a subsequence of $\alpha_n $ (still denoted by $\alpha_n $) such that}
\begin{eqnarray}
\label{c1}
u^{\alpha_n}\!\!\!\!&\to u&\mbox{weakly-star in $\liq$ \pier{and weakly in} $H^1(0,T;H)$},\qquad\\[1mm]
\label{c2}
\mu^{\alpha_n}\!\!\!\!&\to \mu&\mbox{weakly-star in $\liq\pier{{}\cap{}} L^\infty(0,T;V)$}\nonumber\\
&{}&\mbox{and weakly in $H^1(0,T;H)\cap L^2(0,T;W)$,}\\[1mm]
\label{c3}
\rho^{\alpha_n}\!\!\!\!&\to \rho&\mbox{weakly-star in $\liq\pier{{}\cap{}}L^\infty(0,T;V)$}\nonumber\\
&{}&\mbox{and weakly in $H^1(0,T;H)$},\\[1mm]
\label{c4}
\rho_t^{\alpha_n}\!\!\!\!&\to \rho_t&\mbox{weakly in $L^6(Q)$},\\[1mm]
\label{c5}
\gianni{\varphi(\alpha_n)} h'(\rho^{\alpha_n})\!\!\!&\to \xi&\mbox{weakly in $L^6(Q)$}. 
\end{eqnarray}
At first, we note that the set $\uad$ is bounded, convex, and closed (and thus weakly sequentially 
closed) in $H^1(0,T;H)$, whence it follows that
$u\in\uad$.

We now claim that $(\mu,\rho,\xi)$ coincides with the unique solution to the state system \eqref{ss1}--\eqref{ss5}
established in Theorem~2.1. At first, $\mu\in L^2(0,T;W)$ satisfies the boundary condition \eqref{ss4},
and we have that $\,\Delta\mu^{\alpha_n}\to \Delta\mu\,$ weakly in $L^2(Q)$. Moreover,
we can infer from standard compactness results (cf. \cite[Sect.\,8,~Cor.\,4]{Simon}) that, 
for every $q\in [1,6)$, it holds that
\begin{align}
\label{3.2}
&\rho^{\alpha_n}\to \rho \quad \pier{\mbox{strongly in $C^0([0,T];L^q(\oma))$}},\\[1mm]
\label{3.3}
&\mu^{\alpha_n}\to \mu \quad \pier{\mbox{strongly in $C^0([0,T];L^q(\oma))$}},
\end{align} 
\pier{and both pointwise a.\,e. in $Q$, at least for a subsequence of $\alpha_n $.}
Consequently, $\mu$ and $\rho$ fulfill the initial conditions from \eqref{ss5}, and \eqref{reguss3} holds true. 
In addition, we conclude from {\bf (A1)} that, for $\Phi\in\{g,g',F'\}$ and $1\le q<6$,
\beq
\label{3.4}
\Phi(\rho^{\alpha_n}) \to \Phi(\rho) \quad\pier{\mbox{strongly in $C^0([0,T];L^q(\oma))$}}
\eeq
\pier{and pointwise a.\,e. in $Q$.}
From \eqref{c2}--\eqref{3.4}, it readily follows that we have the \gianni{convergence given~by}
\begin{align}
g(\rho^{\alpha_n})\mu_t^{\alpha_n}\to g(\rho)\mu_t,\quad \mu^{\alpha_n}g'(\rho^{\alpha_n})\rho_t^{\alpha_n}
\to \mu g'(\rho)\rho_t,\quad\mu^{\alpha_n}g'(\rho^{\alpha_n})\to\mu g'(\rho),\nonumber\\[1mm]
\mbox{all weakly in $L^1(Q)$, as $n\to\infty$.} 
\end{align}
Moreover, assumption {\bf (A3,vi)} (see, in particular, condition \eqref{B5}) implies that $\cb$ is Lipschitz
continuous on $\lzq$, which entails that
\beq
\label{3.6}
\cb[\rho^{\alpha_n}]\to \cb[\rho] \quad\mbox{strongly in $\lzq$.}
\eeq

In conclusion, the triple $(\mu,\rho,\xi)$ also fulfills \eqref{ss1} and \eqref{ss2}. It remains to show that 
\eqref{ss3} holds true as well. 
Once this \gianni{is} shown, we can infer 
that $(\mu,\rho) = \cs_0(u)$, i.\,e., that the pair $\,((\mu,\rho,\xi),u)\,$ is admissible for $({\bf P}_0)$. 

\vspace{2mm}
To this end, we recall that $\,h\,$ is convex and bounded in $[0,1]$
\gianni{and that $\varphi$ is nonnegative}. 
We thus have, for every $n\in\nz$,
\begin{eqnarray}
\label{3.7}
&& \gianni{\varphi(\alpha_n)\min_{[0,1]}h \cdot \mathop{\rm meas} Q
+ \texinto\varphi(\alpha_n)\,h'(\rho^{\alpha_n})\, (z-\rho^{\alpha_n})\dx\dt}\nonumber\\[1mm]
&&\gianni{{}\leq{}}\texinto\varphi(\alpha_n)\,h(\rho^{\alpha_n})\dx\dt\,
+\,\texinto\varphi(\alpha_n)\,h'(\rho^{\alpha_n})\, (z-\rho^{\alpha_n})\dx\dt\nonumber\\[1mm]
&&\leq\, \texinto\varphi(\alpha_n)\,h(z)\dx\dt\nonumber\\[2mm]
&&\mbox{for all \,}
z\in {\cal K}=\{v\in {L^2(Q)}:0\le v\leq 1\text{ a.\,e. in }Q\}\,. 
\end{eqnarray}
\gianni{Now, we notice that 
the first term and the last one of this chain tend to zero as $n\to\infty$ by~(\ref{defphi})}.
Hence, by \gianni{ignoring the middle line of \eqref{3.7} and} invoking (\ref{c5}) and \eqref{3.2},
passage to the limit as $n\to\infty$ yields
\begin{equation*}
\texinto\xi\,( \rho-z)\dx\dt\,\geq 0\quad\forall \,z\in {\mathcal{K}} . 
\end{equation*}
This entails that $\xi$ is an element of the subdifferential of the extension $\mathcal{I} $ of $ I_{[0,1]}$ to $L^2(Q)$, which means that $\xi \in \partial \mathcal{I}(\rho)$ or, equivalently (cf.~\cite[Ex.~2.3.3., p.~25]{Brezis}),  
that $\xi\in\partial I_{[0,1]}(\rho)$ a.\,e. in $Q$. 

\vspace{2mm}
It remains to show that $((\mu,\rho,\xi),u)$ is in fact optimal for $({\bf P}_0)$.
To this end, let $v\in\uad$ be arbitrary. In view of the convergence properties \gianni{\eqref{c1}--\eqref{3.3}},
and using the weak sequential lower semicontinuity properties of the cost functional, we obtain 
\begin{eqnarray}
\label{3.8}
&& J((\mu,\rho),u)\,=\, J({\cal S}_{0}(u),u)\,\le\,
\liminf_{n\to\infty}\,J(\cs_{\alpha_n}(u^{\alpha_n}),u^{\alpha_n})\nonumber\\[1mm]
&&\leq\,\liminf_{n\to\infty}\, J({\cal S}_{\alpha_n}(v),v) 
\,=\,\lim_{n\to\infty} J({\cal S}_{\alpha_n}(v),v)\,=\,
J({\mathcal S}_{0}(v),v),
\end{eqnarray}  
where for the last equality the continuity in $\lzq\times\lzq$ of the cost functional with respect to 
$(\mu,\rho)$ was used \pier{(see however the next statement)}. With this, the assertion is completely proved.\qed

\vspace{5mm}
{\sc Corollary 3.2.}\quad\,{\em Let the general assumptions} {\bf (A1)}--{\bf (A4)} {\em satisfied, and
let sequences {$\, \{\alpha_n\}\subset (0,1]\,$ and $\,\{u^{\alpha_n}\}\subset {\cal U}$}  be given such that,
as $n\to\infty$, $\,\alpha_n\searrow 0\,$ and $\,u^{\alpha_n}\to u\,$ weakly-star in $\,{\cal X}$. 
Then, with 
$(\mu^{\alpha_n},\rho^{\alpha_n})=\cs_{\alpha_n}(u^{\alpha_n})$, $n\in\nz$, and $(\mu,\rho)=\cs_0(u)$, 
\gianni{\eqref{c2}--\eqref{3.6}} hold true, where $\xi\in L^6(Q)$ satisfies $\xi\in\partial I_{[0,1]}(\rho)$ 
almost everywhere in $Q$. 
Moreover,
we have}
\begin{eqnarray}
\label{3.9}
\lim_{n\to\infty}  J({\cal S}_{\alpha_n}(v),v)\,=\, J({\cal S}_0(v),v) \quad\forall \,v\in\uad\,.
\end{eqnarray}

\vspace{3mm}
{\sc Proof:}\quad\,By the same arguments as in the first part
of the proof of Theorem~3.1, we can conclude that \gianni{\eqref{c2}-\eqref{3.6}}
hold true at least for some subsequence. 
But, as we have \gianni{just seen},
the limit is given by the unique solution triple to the state system 
\eqref{ss1}--\eqref{ss5}.  
Hence, the limit is the same for all convergent subsequences, and thus 
\eqref{c2}--\eqref{3.6} \pier{are} true for the entire sequence, as claimed. 

Now, let $v\in\uad$ be arbitrary.
Then, owing to \eqref{3.2}, \eqref{3.3}, ${\cal S}_{\alpha_n}(v)\,$ converges strongly to $\,{\cal S}_0(v)\,$ in 
$L^2(Q)\times L^2(Q)$. The validity of \eqref{3.9} is then a consequence of the fact that $J$ is continuous
in $L^2(Q)\times L^2(Q)$ with respect \gianni{to~$(\mu,\rho)$}.\qed

\vspace{7mm}
Theorem~3.1 does not yield any information on whether every solution to the optimal control problem $({\bf{P}}_{0})$ can be approximated by a sequence of solutions to the problems $({\bf {P}}_{\alpha})$. 
As already announced in the introduction, we are not able to prove such a general `global' result. Instead, we 
give a `local' answer for every individual optimizer of $({\bf {P}}_{0})$. For this purpose,
we employ a trick due to Barbu~\cite{Barbu}. Now let $\bar u\in\uad$
be an arbitrary optimal control for $({\bf P}_0)$, and let $(\bar \mu,\bar \rho,\bar\xi)$
be the associated solution triple to the state system (\ref{ss1})--(\ref{ss5}) in the sense of 
Theorem 2.1. In particular, $\,(\bar \mu,\bar \rho)={\cal S}_0 (\bar u)$. We associate with this 
optimal control the `adapted cost functional'
\begin{equation}
\label{acost}
\widetilde J((\mu,\rho),u):= J((\mu,\rho\pier{)},u)\,+\,\frac{1}{2}\,\|u-\bar{u}\|^2_{L^2(Q)}
\end{equation}
and a corresponding `adapted optimal control problem'
\begin{eqnarray*}
(\widetilde{\bf {P}}_{\alpha})\quad\mbox{Minimize }\,\, \widetilde J((\mula,\rla),u)\quad
\mbox{subject to $u\in\uad$ and \eqref{sa1}--\eqref{sa4}.}
\end{eqnarray*}

\vspace{3mm}
With a standard direct argument that needs no repetition here, we can show the following 
result.

\vspace{5mm}
{\sc Lemma 3.3.}\quad\,{\em Suppose that the assumptions} {\bf (A1)}--{\bf (A4)} 
{\em are fulfilled. Then the optimal control problem $(\widetilde{\bf P}_\alpha)$
has for every $\alpha\in (0,1]$ a solution.}
  
\vspace{5mm}
We are now in the position to give a partial answer to the question raised above. We have the following result.

\vspace{5mm}
{\sc Theorem~3.4.}\,\quad{\em Let the general assumptions} {\bf (A1)}--{\bf (A4)} {\em be fulfilled,
suppose that $\bar u\in \uad$ is an arbitrary optimal control of} $({\bf P}_{0})$ {\em with associated state
triple $(\bar\mu,\bar\rho,\bar\xi)$, and let $\,\{\alpha_n\}\subset (0,1]$ be any sequence such that
$\,\alpha_n\searrow 0\,$ as $\,n\to\infty$. Then there exist a subsequence $\{\alpha_{n_k}\}$
of $\{\alpha_n\}$ and, for every $k\in\nz$, an optimal control
 $\,u^{\alpha_{n_k}}\in \uad\,$ of the adapted problem $(\widetilde{\bf P}_{\alpha_{n_k}})$
with associated state $(\mu^{\alpha_{n_k}},\rho^{\alpha_{n_k}})$ such that, as $k\to\infty$,}
\gianni{{\em
\begin{eqnarray}
\label{3.10}
&&u^{\alpha_{n_k}}\to \bar u\quad\mbox{strongly in }\,L^2(Q),\\[1mm]
\label{3.11}
&&\mbox{\eqref{c2}--\eqref{3.6} hold true, where $(\mu,\rho,\xi)$ is
replaced} \nonumber\\
&&\mbox{by $(\bar\mu,\bar\rho,\bar\xi)$ and the index $n$ is replaced by $n_k$,}\\[1mm]
\label{3.12}
&&\widetilde J((\mu^{\alpha_{n_k}},\rho^{\alpha_{n_k}}),u^{\alpha_{n_k}})\to  J((\bar\mu,\bar\rho),\bar u)\,.
\end{eqnarray}
}}

{\sc Proof:} \quad\,Let $\alpha_n \searrow 0$ as $n\to\infty$. 
For any $ n\in\nz$, we pick an optimal control 
$u^{\alpha_n} \in \uad\,$ for the adapted problem $(\widetilde{\bf P}_{\alpha_n}\gianni)$ and denote by 
$(\mu^{\alpha_n},\rho^{\alpha_n})={\cal S}_{\alpha_n}(u^{\alpha_n})$ the associated solution to  
problem (\ref{sa1})--(\ref{sa4}), where the right-hand side in \eqref{sa1} is given by $u^{\alpha_n}$. 
Then, by virtue of Theorems 2.2 \gianni{and~2.3}, (\ref{regusa1})--(\ref{sabounds1}) and 
\eqref{sabounds2}--\eqref{sabounds3} are satisfied. Arguing as in the proof of Corollary~3.2,
and using the fact that $\uad$ is a bounded subset of~${\cal X}$ \gianni{(see~\eqref{defX})}, 
we can infer that there is a subsequence $\{n_k\}_{k\in\nz}$
of $\nz$ such that
\begin{equation}
\label{3.13neverquoted}
u^{\alpha_{n_k}}\to u\quad\mbox{weakly-star in }\,{\cal X}
\quad\mbox{as }\,k\to\infty,
\end{equation}
with some $u\in\uad$, and such that \eqref{3.11} holds true with $(\mu,\rho):=\cs_0(u)$
and a suitable $\xi\in\partial I_{[0,1]}(\rho)$. 
In particular, the pair $((\mu,\rho,\xi),u)$
is admissible \gianni{for~$({\bf P}_0)$}.

\vspace{2mm}
We now aim to prove that $u=\bar u$. Once this \gianni{is} shown, the uniqueness result of Theorem~2.1 yields 
that also $\,(\mu,\rho,\xi)=(\bar\mu,\bar\rho,\bar\xi)$, which then implies that (\ref{3.11}) 
holds true. 
Indeed, we have, owing to the weak sequential lower semicontinuity of $\widetilde J$, 
and in view of the optimality property of $\,((\bar\mu,\bar\rho),\bar u)$ for problem $({\bf P}_0)$,
\begin{eqnarray}
\label{3.13}
&&\liminf_{k\to\infty}\, \widetilde  J((\mu^{\alpha_{n_k}},\rho^{\alpha_{n_k}}), u^{\alpha_{n_k}})
\ge \, J((\mu,\rho),u)\,+\,\frac{1}{2}\,
\|u -\bar{u}\|^2_{L^2(Q)}\nonumber\\[1mm]
&&\geq \, J((\bar\mu,\bar\rho),\bar u)\,+\,\frac{1}{2}\,\|u-\bar{u}\|^2_{L^2(Q)}\,.
\end{eqnarray}
On the other hand, the optimality property of  $\,((\mu^{\alpha_{n_k}},\rho^{\alpha_{n_k}}),
u^{\alpha_{n_k}})\,$ for problem $(\widetilde {\bf P}_{\alpha_{n_k}})$ yields that
for any $k\in\nz$ we have
\begin{equation}
\label{3.14}
\widetilde J((\mu^{\alpha_{n_k}},\rho^{\alpha_{n_k}}),
u^{\alpha_{n_k}})\, =\,
\widetilde J\left({\cal S}_{\alpha_{n_k}}(u^{\alpha_{n_k}}),
u^{\alpha_{n_k}}\right)\,\le\,\widetilde J\left({\cal S}_{\alpha_{n_k}}(\bar u),\bar u\right)\,,
\end{equation}
whence, taking the limit superior as $k\to\infty$ on both sides and invoking \eqref{3.9} in
Corollary~3.2,
\begin{align}
\label{3.15}
&\limsup_{k\to\infty}\,\widetilde J((\mu^{\alpha_{n_k}},\rho^{\alpha_{n_k}}),
u^{\alpha_{n_k}})\,\le\,\widetilde J({\cal S}_0(\bar u),\bar u)\nonumber\\[1mm]
&=
\pier{\widetilde J((\bar\mu,\bar\rho),\bar u)={}}J((\bar\mu,\bar\rho),\bar u)\,.
\end{align}
Combining (\ref{3.13}) with (\ref{3.15}), we have thus shown that 
$\,\frac{1}{2}\,\|u-\bar u\|^2_{L^2(Q)}=0$\,,
so that in fact $\,u=\bar u\,$  and thus also $\,(\mu,\rho,\xi)=(\bar\mu,\bar\rho,\bar\xi)$.
Moreover, (\ref{3.13}) and (\ref{3.15}) also imply that
\begin{eqnarray}
&&J((\bar\mu,\bar\rho),\bar u) \, =\,\widetilde J((\bar\mu,\bar\rho),\bar u)
\,=\,\liminf_{k\to\infty}\, \widetilde J((\mu^{\alpha_{n_k}},\rho^{\alpha_{n_k}}),
u^{\alpha_{n_k}})\nonumber\\[1mm]
&&\,=\,\limsup_{k\to\infty}\, \widetilde J((\mu^{\alpha_{n_k}},\rho^{\alpha_{n_k}}),u^{\alpha_{n_k}})\,
=\,\lim_{k\to\infty}\, \widetilde J((\mu^{\alpha_{n_k}},\rho^{\alpha_{n_k}}),
u^{\alpha_{n_k}})\,,\qquad
\end{eqnarray}                                     
which proves \eqref{3.12} and, at the same time, also (\ref{3.10}). The assertion is thus
completely proved.\qed

\section{The optimality system}
\setcounter{equation}{0}
In this section, we aim to establish first-order necessary optimality conditions for the optimal control problem $({\mathcal{P}}_{0})$.  
This will be achieved by a passage to the limit as $\alpha\searrow 0$ 
\gianni{in the first-order necessary optimality conditions for the adapted optimal control problems $(\widetilde{\bf{P}}_{\alpha})$
that can be derived by arguing as in~\cite{CGS4} mith minor changes}. 
\gianni{This procedure will yield certain generalized first-order necessary conditions of optimality in the limit}. 
In this entire section, we assume that $\bar u\in\uad$ is a fixed optimal control 
of problem $({\bf P}_0)$ and that $(\bar\mu,\bar\rho,\bar\xi)$ is the associated solution
to the state system \eqref{ss1}--\eqref{ss5} established in Theorem 2.1, that is, we have $(\bar\mu,\bar\rho)
=\cs_0(\bar u)$ and $\xi\in\partial I_{[0,1]}(\rho)$ almost everywhere in $Q$. 

We begin our analysis by formulating for arbitrary $\alpha \in (0,1]$ the adjoint state system for the adapted
control problem $(\widetilde{\bf{P}}_{\alpha})$ corresponding to $\bar u$.
We assume that $\,\bar u^\alpha\in\uad\,$ is an arbitrary optimal control for $(\widetilde{\bf P}_\alpha)$
and that $\,(\bar\mu^\alpha,\bar\rho^\alpha)=\cs_\alpha(\bar u^\alpha)$ is the corresponding solution to 
the associated state system (\ref{sa1})--(\ref{sa4}), which then enjoys the regularity properties 
\eqref{regusa1}, \eqref{regusa2} and fulfills the boundedness conditions \eqref{sabounds1}, \eqref{sabounds2}
and \eqref{sabounds3}.  It then follows (see \cite[Eqs.~(4.3)--(4.6)]{CGS4}) that the corresponding adjoint
system has the form
\begin{align}
\label{adj1}
&-(1+2g(\bar\rho^\alpha))\,p_t^\alpha - g'(\bar\rho^\alpha)\,\bar\rho_t^\alpha\,\pla-\Delta\pla -
g'(\bar\rho^\alpha)\,\qla\nonumber\\[1mm]
&\quad=\,\gian{\beta_2(\bar\mu^\alpha-\mu_Q)}
\quad\mbox{in \,$Q$,}
\\
\label{adj2}
&-q^\alpha_t + \vp\,h''(\rho^\alpha)\,q^\alpha +F''(\bar\rho^\alpha)\,\qla-\bar\mu^\alpha\,g''(\bar\rho^\alpha)\,\qla
\nonumber\\[1mm] 
&\quad=\,- g'(\bar\rho^\alpha)(\bar\mu_t^\alpha\,\pla-\bar\mu^\alpha\,p_t^\alpha) - D\cb[\bar\rho^\alpha]^*(\qla)
+ \gian{\beta_1(\bar\rho^\alpha-\rho_Q)}
\quad\mbox{in \,$Q$,}\\[1mm]
\label{adj3}
&\qquad\dn\pla=0\quad\mbox{on \,$\Sigma$}, \\[1mm]
\label{adj4}
&\qquad \pla(T)=\qla(T)=0 \quad\mbox{in \,$\oma$},   
\end{align}
where $D\cb[\bar\rho^\alpha]^*\in {\cal L}(\lzq,\lzq)$ denotes the adjoint operator associated with the
operator $D\cb[\bar\rho^\alpha]\in {\cal L}(\lzq,\lzq)$, which is defined 
by the identity
\beq
\label{adjoint}
\texinto D\cb[\bar\rho^\alpha]^*(v)\,w\dx\dt\,=\texinto v\,D\cb[\bar\rho^\alpha](w)\dx\dt\quad
\forall\,v,w\in\lzq.
\eeq

According to \cite[Thm.~4.2]{CGS4}, the system (\ref{adj1})--(\ref{adj3})
has for every $\alpha \in (0,1]$ a unique solution \gianni{pair} $(p^\alpha,q^\alpha)$ such that
\begin{equation}
\label{reguadj}
p^\alpha\in H^1(0,T;H)\cap L^\infty(0,T;V)\cap L^2(0,T;W),\quad
q^\alpha \in H^1(0,T;H).
\end{equation}

Moreover, as in \cite[Cor.~4.3]{CGS4} it follows that the following variational 
inequality is satisfied:
\begin{equation}
\label{vug1}
\texinto \bigl(p^\alpha + \beta_3\,\bar u^\alpha + (\bar u^\alpha-\bar u)\bigr)
(v-\bar u^\alpha)\dx\dt\,\ge\,0 \quad
\forall\,v \in \uad\,.
\end{equation}
 
\vspace{3mm}
We now prove an a priori estimate that will be fundamental for the derivation of the optimality conditions for 
$({\bf P}_0)$. To this end, we set, for all $\alpha\in (0,1]$,
\beq
\label{lamal}
\lambda^\alpha:=\vp \,h''(\bar\rho^\alpha)\,q^\alpha,
\eeq
and we introduce the function space
\begin{equation}
\label{defY}
Y\,:=\,\{v\in H^1(0,T;H): v(0)=0\},
\end{equation}
which is a Hilbert space when equipped with the standard inner product of $H^1(0,T;H)$. 
\gianni{As $H^1(0,T;H)$} is continuously embedded in $C^0([0,T];H)$, the initial condition
encoded in the definition of $Y$ is meaningful. We denote by $Y'$ the dual space of $Y$   
and by  $\langle \, \cdot\, ,  \, \cdot\, \rangle_Y$
the {duality} pairing between $Y'$ and~$Y$. Observe that $Y\subset L^2(Q)\subset Y'$
with dense and continuous embeddings, where it is understood that for all
$w\in Y$ and $v\in L^2(Q)$ it holds that
\beq
\label{duality}
\langle v,w\rangle_Y \,=\, \texinto v\,w \dx\dt\,.
\eeq

We have the following result.

\vspace{5mm}
{\sc Lemma 4.1.} \quad {\em Suppose that the conditions} {\bf (A1)}--{\bf (A4)} {\em are fulfilled.
Then there is a constant $K_3^*>0$, which depends only on the data of problem $({\bf P}_0)$, such that,
for every $\alpha\in (0,1]$,}
\begin{align}
\label{estidual}
&\left\|p^\alpha\right\|_{H^1(0,T;H)\cap L^\infty(0,T;V)\cap L^2(0,T;H^2(\oma))}
\,+\,\left\|q^\alpha\right\|_{L^\infty(0,T;H)}\nonumber\\[1mm]
&+\,\left\|\pt q^\alpha\right\|_{Y'}\,+\,
\left\|\lambda^\alpha\right\|_{Y'}\,\le\,K_3^*\,.
\end{align}

\vspace{3mm}
{\sc Proof:} \quad\,In the following, $C_i$, $i\in\nz$, denote positive constants, which are independent
of $\alpha\in (0,1]$. For the sake of a better readability, we will in the following steps omit the 
superscript $\alpha$, writing it only at the final estimate of each step. We also will make repeated
use of the global bounds \eqref{sabounds1}, \eqref{sabounds2}, \eqref{sabounds3} without further reference.

\vspace{5mm}
\underline{\sc Step 1:} \quad We first add $\,p\,$ on both sides of \eqref{adj1}, multiply the result by $-p_t$, 
and integrate over~$Q^t$ \gianni{(recall \eqref{defQt})}, 
where $0\le t<T$. Using the fact that $g$ is nonnegative, we then obtain that
\beq
\label{step11}
\itt |p_t|^2\dx\ds\,+\,\frac 12\,\|p(t)\|^2_V\,\le\,\sum_{j=1}^4 I_j,
\eeq
where the quatities $I_j$, $1\le j\le 4$, will be specified and estimated below. At first, we obtain from
Young's inequality the estimates
\begin{align}
\label{step12}
&I_1:=\,-\itt p\,p_t\dx\ds\,\le\,\frac 18\itt |p_t|^2\dx\ds\,+\,C_1\itt\!\pier{|p|}^2\dx\ds\,,\\[2mm]
\label{step13}
&I_3:=\,-\itt g'(\bar\rho)\,q\,p_t\dx\ds\,\nonumber\\[1mm]
&\hspace*{4.3mm}\le\,\frac 18\itt|p_t|^2\dx\ds\,+\,C_2\itt|q|^2\dx\ds\,,
\\[2mm]
\label{step14}
&I_4:=\,-\itt\gian{\beta_2(\bar\mu-\mu_Q)}\,p_t\dx\ds\,
\le\,\frac 18 \itt|p_t|^2\dx\ds\,+\,C_3\,.\qquad
\end{align} 
Moreover, owing to H\"older's and Young's inequality, and using the continuity of the embedding $V\subset L^6(\oma)$,
\begin{align}
\label{step15}
&I_2:=\,-\itt g'(\bar\rho)\,\bar\rho_t\,p\,p_t\dx\ds\,\le\,C_4\int_t^T\!\|p_t(s)\|_2\,\|\bar\rho_t(s)\|_3\,
\|p(s)\|_6\ds\nonumber\\[1mm]
&\hspace*{4mm}\, \le \,\frac 18\itt |p_t|^2\dx\ds\,+\,C_5\int_t^T\!\|\bar\rho_t(s)\|_3^2\,\|p(s)\|_V^2\ds\,.
\end{align}

Combining the estimates \eqref{step11}--\eqref{step15}, we have thus shown that
\begin{align}
\label{step16}
& \itt|p_t^\alpha|^2\dx\ds\,+\,\|p^\alpha(t)\|_V^2\nonumber\\[1mm]
&\le\,C_6\Big(1+\itt \bigl(|p^\alpha|^2+|q^\alpha|^2\bigr)\dx\ds\,+\,\int_t^T\|\bar\rho^\alpha_t(s)\|_3^2\,
\|p^\alpha(s)\|_V^2\ds
\Big)\,.
\end{align}

\vspace{5mm}
\underline{\sc Step 2:} \quad We now multipy \eqref{adj2} by $q$ and integrate over $Q^t$, where
$0\le t<T$. We obtain that 
\beq\label{step21}
\frac 12\,\|q(t)\|_H^2\,+\itt (\vp\,h''(\bar\rho)-\mu\,g''(\bar\rho))\,|q|^2\dx\ds\,\le\,\sum_{j=1}^5 J_j,
\eeq
where, since $\gianni{\vp}h''(\bar\rho)\ge 0$ and $\mu\,g''(\bar\rho)\le 0$, the integral on the left-hand side
is nonnegative and  the quantities $J_j$, $1\le j\le 5$, will be defined \pier{and} estimated below. We have
\begin{equation}
\label{step22}
J_1:=\,-\itt\!F''(\bar\rho)|q|^2\dx\ds\,\le\,C_7\itt\!|q|^2\dx\ds,
\end{equation}
as well as, using Young's inequality,
\begin{align}
\label{step23}
&J_3:=\,\itt\! g'(\bar\rho)\,\mu\,p_t\,q\dx\ds\nonumber\\[1mm]
&\hspace*{5mm}\le\,\frac 14\itt\!|p_t|^2\dx\ds\,+\,C_8\itt\!|q|^2\dx\ds,\\[2mm]
\label{step24}
&J_5:=\,\itt\!\gian{\beta_1(\bar\rho-\rho_Q)}\,q\dx\ds\,
\le\,\itt\!|q|^2\dx\ds\,+\,C_9,\\[2mm]
\label{step25}
&J_4:=\,-\itt D\cb[\bar\rho]^*(q)\,q\dx\ds\,=\,-\itt\!D\cb[\bar\rho](q)\,q\dx\ds\nonumber\\[1mm]
&\hspace*{5mm}\le\,C_{10}\itt\!|q|^2\dx\ds,
\end{align}
where in the last estimate we have used \eqref{adjoint} and \eqref{B5}. Finally, we 
\pier{take advantage} of
the continuity of the embedding $H^2(\oma)\subset \lio$ and of H\"older's and Young's inequality
to conclude that, for every $\gamma>0$ (to be chosen later),
\begin{align} 
\label{step26}
&J_2:=\,-\itt\! g'(\bar\rho)\,\bar\mu_t\,p\,q\dx\ds\,\le\,C_{11}
\int_t^T\! \|\bar\mu_t(s)\|_2\,\|p(s)\|_\infty\,\|q(s)\|_2\ds\nonumber\\[1mm]
&\hspace*{5mm}\le\,\gamma\int_t^T\!\|p(s)\|_{H^2(\Omega)}^2\ds\,+\,\frac {C_{12}}\gamma
\int_t^T\!\|\bar\mu_t(s)\|_2^2\,\|q(s)\|_2^2\ds\,.
\end{align}

In view of the boundary condition \eqref{adj3}, we can infer from classical elliptic estimates that
\begin{equation}
\label{step27}
\|p(s)\|_{H^2(\oma)}\,\le\,C_\oma \,(\|p(s)\|_H\,+\,\|\Delta p(s)\|_H)\quad\mbox{for a.\,e. }\,s\in (0,T),
\end{equation}
where the constant $C_\oma>0$ depends only on the domain $\oma$.

Combining all the estimates \eqref{step21}--\eqref{step27}, we have therefore shown that
\begin{align}
\label{step28}
&\|q^\alpha(t)\|_H^2\,\le\,C_{13}\,+\,\frac 12 \itt\!|p_t^\alpha|^2\dx\ds
\,+\,\gamma\,C_{14}\itt\!\bigl(|p^\alpha|^2+|\Delta p^\alpha|^2\bigr)\dx\ds\nonumber\\[1mm]
&\hspace*{21.5mm} +\,C_{15}\left(1+\gamma^{-1}\right)\itt\!\bigl(\gianni{\|\bar\mu_t^\alpha(s)\|_2^2+1}\bigr)\,
|q^\alpha|^2\dx\ds\,.
\end{align}

\vspace{5mm}
\underline{\sc Step 3:} \quad We now use \eqref{adj1} to estimate $\Delta p^\alpha$ directly. Indeed, we have
\begin{align}
\label{step31}
&\itt|\Delta p^\alpha|^2\dx\,ds\,=\,\|\Delta p^\alpha\|_{L^2(Q^t)}^2 \nonumber\\[1mm]
&\le\,C_{16}\,+\,C_{17}\left(\|p_t^\alpha\|_{L^2(Q^t)}^2 \,+\,\|q^\alpha\|_{L^2(Q^t)}^2
\,+\,\|\bar\rho^\alpha_t\,p^\alpha\|_{L^2(Q^t)}^2\right).
\end{align}
Now, owing to \eqref{embed2}, it holds that
\beq
\label{step32}
\| \bar\rho^\alpha_t\,p^\alpha\|_{L^2(Q^t)}^2\,\le \,C_{18}\int_t^T \|\bar\rho^\alpha_t(s)\|_3^2\,
\|p^\alpha(s)\|_V^2\ds\,.
\eeq
In conclusion, we have shown that
\begin{align}
\label{step33}
&\itt\!|\Delta p^\alpha|^2\dx\ds\,\le\,C_{16}\,+\,C_{17}\itt\bigl(|p_t^\alpha|^2+|q^\alpha|^2\bigr)\dx\ds
\nonumber\\[1mm]
&\hspace*{39mm} +\,C_{19}\int_t^T \|\bar\rho^\alpha_t(s)\|_3^2\,\|p^\alpha(s)\|_V^2\ds\,.
\end{align}

\vspace{4mm}
\underline{\sc Step 4:} \quad At this point, we are in a position to combine the estimates \eqref{step16},
\eqref{step28} and \eqref{step33}. Choosing $\gamma>0$ appropriately small (such that 
$0<\gamma\,C_{14}\,C_{17}<\frac 12$), and noting that the global estimates \eqref{sabounds3} imply that
the mappings $\,\,s\mapsto\|\bar\rho^\alpha_t(s)\|_3^2\,\,$ and $\,\,\|\bar\mu^\alpha_t(s)\|_2^2\,\,$
are bounded in $L^1(0,T)$ independently of $\alpha\in (0,1]$, we conclude from Gronwall's lemma (taken
backward in time) the estimate
\beq
\label{step41}
\left\|p^\alpha\right\|_{H^1(0,T;H)\cap L^\infty(0,T;V)\cap L^2(0,T;H^2(\oma))}\,+\,
\left\|q^\alpha\right\|_{L^\infty(0,T;H)}\,\le\,C_{20}.
\eeq 

Now assume that $v\in Y$ is arbitrary. As $q^\alpha(T)=v(0)=0$, integration by parts with respect
to time (which is permitted since $q^\alpha,v\in H^1(0,T;H)$) yields 
\begin{align}
&\left|\langle q^\alpha_t,v\rangle_Y\right|\,=\,\Big|\texinto q^\alpha_t\,v\dx\dt\Big|\,=\,
\Big|-\texinto q^\alpha\,v_t\dx\dt \Big|\nonumber\\[2mm]
&\hspace*{18mm}\le\,C_{21}\left\|q^\alpha\right\|_{L^\infty(0,T;H)}\,\|v\|_{H^1(0,T;H)},
\end{align}  
and it follows from \eqref{step41} that
\beq
\label{step42}
\left\|q^\alpha_t\right\|_{Y'}\,\le\,C_{22}.
\eeq
It remains to show the bound for $\lambda^\alpha$ \gianni{(see \eqref{lamal})}. 
To this end, notice that \gianni{\eqref{adj2} yields}
\begin{align}
&\lambda^\alpha\,=\,q^\alpha_t \,-\,g'(\bar\rho^\alpha)\,\bar\mu_t^\alpha\,p^\alpha
-F''(\bar\rho^\alpha)\,q^\alpha\,+\,\bar\mu^\alpha\,g''(\bar\rho^\alpha)\,q^\alpha
\,+\,\bar\mu^\alpha\,g'(\bar\rho^\alpha)\,p^\alpha_t\nonumber\\[1mm]
&\hspace*{9.5mm} -\,D\cb[\bar\rho^\alpha]^*(q^\alpha)\,+\,\gian{\beta_1(\bar\rho^\alpha-\rho_Q)},
\end{align}
where it easily follows from the estimates \eqref{sabounds2}, \eqref{sabounds3}
and \eqref{step41} that the last five summands on the right-hand side are bounded
in $L^2(Q)$ (and thus \gianni{in~$Y'$}) independently of $\alpha\in (0,1]$. 
To estimate the remaining term, let $v\in Y$ be arbitrary
\gianni{and observe that \eqref{regusa1}, \eqref{reguadj}, 
the continuous embedding $V\subset L^6(\Omega)$ and H\"older's inequality imply
$$
  \bar\mu_t^\alpha\,p^\alpha
  \in L^\infty(0,T;L^{3/2}(\Omega)) \cap L^1(0,T;L^3(\Omega))
  \subset L^2(0,T;L^2(\Omega)) ,
$$
whence also}
$g'(\bar \rho^\alpha)\,\bar\mu_t^\alpha\,p^\alpha\in L^2(Q)$.
\gianni{Therefore}, using \eqref{duality}, \eqref{sabounds3} and \eqref{step41}, as well as the continuity of the embeddings
$H^2(\oma)\subset \lio$ and $H^1(0,T;H)\subset C^0([0,T];H)$,
\gianni{we have}
\begin{align*}
&\left|\langle -g'(\bar \rho^\alpha)\,\bar\mu_t^\alpha\,p^\alpha,v\rangle_{Y}\right|
\,=\,\Big|\texinto g'(\bar \rho^\alpha)\,\bar\mu_t^\alpha\,p^\alpha\,v\dx\dt\Big|\nonumber\\[1mm]
&\le\,C_{23}\int_0^T\!\|\bar\mu^\alpha_t(t)\|_2\,\|p^\alpha(t)\|_\infty\,\|v(t)\|_2\dt\nonumber\\[1mm]
&\le\,C_{24}\,\max_{0\le t\le T}\,\|v(t)\|_H\,\left\|\bar\mu_t^\alpha\right\|_{L^2(Q)}\,
\left\|p^\alpha\right\|_{L^2(0,T;H^2(\oma))}\nonumber\\[1mm]
&\le\,C_{25}\,\|v\|_{H^1(0,T;H)}\,.
\end{align*}
In consequence, we have that
\beq
\label{step43}
\left\|\lambda^\alpha\right\|_{Y'}\,\le\,C_{26},
\eeq 
which concludes the proof of the assertion.\qed

\vspace{3mm}
We draw some consequences from the previously established results. Assume that $\{\alpha_n\}
\subset (0,1]$ satisfies  $\alpha_n\searrow 0$ as $n\to\infty$. Then, thanks to Theorem~3.4,
there is a subsequence, without loss of generality $\{\alpha_n\}$ itself, such that,  
for any $n\in\nz$, we can find an optimal control $\bar u^{\alpha_n}\in\uad$ for
$(\widetilde{\bf{P}}_{\alpha_n})$ and an associated state $(\bar\mu^{\alpha_n},\bar\rho^{\alpha_n})$ 
such that \gianni{the convergence given by (\ref{3.10}) and (\ref{3.12}) 
(where one reads $\alpha_n$ in place of~$\alpha_{n_k}$) holds true}. 
As in the proof of Theorem~3.1, 
we may without loss of generality assume
that, for $1\le q<6$ \gianni{and $\Phi\in\{F'', g, g',g''\}$}, 
\begin{align}
\label{eq:4.34}
&\Phi(\bar\rho^{\alpha_n})\to \Phi(\bar\rho)\,\quad\mbox{strongly in }\,C^0([0,T];L^q(\oma)).
\end{align}
Also, by virtue of Lemma 4.1, we may without loss of generality assume
that the corresponding adjoint state variables $(p^{\alpha_n},
q^{\alpha_n})$ satisfy 
\begin{align}
\label{eq:4.35}
&p^{\alpha_n}\to p\quad\mbox{weakly-star in $H^1(0,T;H)\cap L^\infty(0,T;\gianni V)\cap L^2(0,T;W)$},\\[1mm]
\label{eq:4.36}
&q^{\alpha_n}\to q \quad\mbox{weakly-star in $L^\infty(0,T;H)$},\\[1mm]
\label{eq:4.38} 
&\lambda^{\alpha_n}\to \lambda \quad\mbox{weakly in $Y'$ },
\end{align}
for \gianni{a suitable triple $\,(p,q,\lambda)$}. 
Therefore, passing to the limit as
$n\to\infty$ in the variational inequality (\ref{vug1}), written for $\alpha_n$, $n\in\nz$, we obtain
that $p$ satisfies
\begin{equation}
\label{vug2}
\texinto\!(p\,+\,\beta_3\,{\bar u})\,(v-{\bar u})\dx\dt\,\ge\,0 \quad
\forall\,v\in\uad.
\end{equation}

\vspace*{2mm}
Next, we will show that in the limit as $n\to \infty$ a limiting adjoint system for $({\cal P}_0)$
is satisfied. To this end, we note that it is not difficult to check that
\begin{align}
\label{eq:4.40}
&g(\bar\rho^{\alpha_n})\,p_t^{\alpha_n}\to g(\bar\rho)\,p_t,\quad g'(\bar\rho^{\alpha_n})\,\bar\rho_t^{\alpha_n}\,
p^{\alpha_n}\to g'(\bar\rho)\,\bar\rho_t\,p,\nonumber\\[1mm]
&g'(\bar\rho^{\alpha_n})\,q^{\alpha_n}\to g'(\bar\rho)\,q,\quad\mbox{all weakly in $L^1(Q)$}. 
\end{align} 
Hence, passage to the limit as $n\to\infty$ in \eqref{adj1} (written for $\alpha_n$) yields 
\begin{align}
\label{nadj1}
&-(1+2g(\bar\rho))\,p_t-g'(\bar\rho)\,\bar\rho_t\,p
-\Delta p-g'(\bar\rho)\,q\nonumber\\[1mm]
&\quad =\,\pier{\beta_2(\bar\mu-\mu_Q)}\quad\mbox{a.\,e. in $\,Q$.} 
\end{align}
We also have that \gianni{(see \eqref{eq:4.35})}
\beq
\label{nadj2}
\dn p=0\quad\mbox{a.\,e. on $\Sigma$}, \quad p(T)=0\quad\mbox{a.\,e. in }\,\oma.
\eeq

In order to derive an equation resembling \eqref{adj2}, we multiply \eqref{adj2} (written for~$\alpha_n$)
by an arbitrary element $v$ belonging to the space
$$
Y_0:=\left\{v\in C^1([0,T];C(\overline{\oma})):v(\cdot,0)=0\right\},$$
which is a dense subset of~$Y$. 
Integrating over $Q$ and by parts with respect to~$t$,
we then obtain the equation
\begin{align}
\label{eq:4.43}
&\texinto\lambda^{\alpha_n}\,v\dx\dt\,+\, \texinto q^{\alpha_n}\,v_t\dx\dt\,+\texinto
F''(\bar\rho^{\alpha_n})\,q^{\alpha_n}\,v\dx\dt\nonumber\\[1mm]
&-\texinto\bar\mu^{\alpha_n}\,g''(\bar\rho^{\alpha_n})\,q^{\alpha_n}\,v\dx\dt
\,+\,\texinto g'(\bar\rho^{\alpha_n})\,\bar\mu_t^{\alpha_n}\,p^{\alpha_n}\,v\dx\dt\nonumber\\[1mm]
&=\,\texinto g'(\bar\rho^{\alpha_n})\,\bar\mu^{\alpha_n}\,p^{\alpha_n}_t\,v\dx\dt
\,-\,\texinto D\cb[\bar\rho^{\alpha_n}]^*(q^{\alpha_n})\,v\dx\dt\nonumber\\[1mm]
&\quad+\texinto\gian{\beta_1(\bar\rho^{\alpha_n}-\rho_Q)}\,v\dx\dt\,.
\end{align}

By virtue of the previously established convergence properties of the involved sequences, 
it is not difficult to show that we may pass to the limit as
$n\to\infty$ in all of the summands occurring in~\eqref{eq:4.43} up to the penultimate one. We may leave the details
of these straightforward calculations to the reader. For the penultimate summand, we have
\begin{align*}
&\texinto \left(D\cb[\bar\rho^{\alpha_n}]^*(q^{\alpha_n})-D\cb[\bar\rho]^*(q)\right)v\dx\dt\nonumber\\[1mm]
&=\texinto q^{\alpha_n}\,(D\cb[\bar\rho^{\alpha_n}](v)-D\cb[\bar\rho](v))\dx\dt
\nonumber\\[1mm]
&\quad +\texinto (q^{\alpha_n}-q)\,D\cb[\bar\rho](v)\dx\dt\,.
\end{align*}
While the second summand tends to zero as $n\to\infty$, nothing can be said about the first one: we need an
additional assumption. The following condition is obviously sufficient to guarantee that also the first summand 
approaches zero as $n\to\infty$:

\vspace{5mm}
{\bf (A5)}  \quad The mapping $\,\,v\mapsto D\cb[v]\,\,$ is continuous \pier{from}
$\lzq$ into ${\cal L}(\lzq,\lzq)$. 

\vspace{5mm}
In conclusion, if {\bf (A5)} is valid, then the passage to the limit as $n\to\infty$ results in the 
following identity:
\begin{align}
\label{eq:4.44}
&\langle \lambda,v\rangle_{\gianni Y}\,+\, \texinto q\,v_t\dx\dt\,+\texinto
F''(\bar\rho)\,q\,v\dx\dt\nonumber\\[1mm]
&-\texinto\bar\mu\,g''(\bar\rho)\,q\,v\dx\dt
\,+\,\texinto g'(\bar\rho)\,\bar\mu_t\,p\,v\dx\dt\nonumber\\[1mm]
&=\,\texinto g'(\bar\rho)\,\bar\mu\,p_t\,v\dx\dt
\,-\,\texinto D\cb[\bar\rho]^*(q)\,v\dx\dt\nonumber\\[1mm]
&\quad+\texinto\gian{\beta_1(\bar\rho-\rho_Q)}\,v\dx\dt\,\qquad\forall\,v\in Y_0.
\end{align}

We claim that the variational equation \eqref{eq:4.44} holds in fact true for all $v\in Y$. To see
this, we employ a standard density argument. Indeed, if $v\in Y$ is given, then there is some
sequence $\{v_n\}\subset Y_0$ such that $v_n\to v$ in the norm of $H^1(0,T;H)$. Now observe that
\eqref{eq:4.44} is valid for $v=v_n$. It is now easily checked that we may pass to the limit as
$n\to\infty$ in each term occurring in \eqref{eq:4.44}. As an example, we give the details for the
most difficult term, which is the last one on the left-hand side. We have, denoting by $C_i$, $i\in\nz$,
constants that do not depend on $n$:
\begin{align*}
&\Big|\texinto g'(\bar\rho)\,\bar\mu_t\,p\,(v-v_n)\dx\dt\Big|\,\le\,C_1\int_0^T\!
\|\bar \mu_t(t)\|_2\,\|p(t)\|_\infty\,\|v(t)-v_n(t)\|_2\,dt\\[1mm]
&\le\,C_2\,\max_{0\le t\le T}\,\|v(t)-v_n(t)\|_H\,\|\bar\mu_t\|_{L^2(Q)}\,\|p\|_{L^2(0,T;H^2(\oma))}\\[1mm]
&\le\,C_3\,\|v-v_n\|_{H^1(0,T;H)}\,\,\to 0\quad\mbox{as $\,n\to\infty$}.
\end{align*}

\vspace*{3mm}
Next, we show that the limit $\,\lambda\,$ satisfies some sort
of a complementarity slackness condition. Indeed, we have
\beq\label{eq:4.45}
\liminf_{n\to\infty}\texinto\lambda^{\alpha_n}\,q^{\alpha_n}\dx\dt\,=\,\liminf_{n\to\infty}
\texinto\varphi(\alpha_n)\,h''(\gianni{\bar\rho}^{\alpha_n})
\,|q^{\alpha_n}|^2\dx\dt\,\ge\,0\,.
\eeq
   
Moreover, there is some indication that the limit $\lambda$ 
should somehow be  concentrated
on the set where $\,\bar\rho=0\,$ or $\,\bar\rho=1$ (which, however, we cannot prove rigorously). 
To this end, we test $\lambda^{\alpha_n}$  by the function
$\,\bar\rho^{\alpha_n}(1-\bar\rho^{\alpha_n})\,\phi$,
where $\,\phi\,$ is any smooth test function satisfying $\phi(0)=0$.
\gianni{By recalling that $h''(r)=\frac 1{r(1-r)}$, we} then obtain that
\begin{align}
\label{eq:4.46}
&\lim_{n\to\infty}\texinto\lambda^{\alpha_n}\,\bar\rho^{\alpha_n}
(1-\bar\rho^{\alpha_n})\,\phi\dx\dt
\,=\,\lim_{n\to\infty}\texinto \varphi(\alpha_n)\,q^{\alpha_n}\,\phi\dx\dt\,=\,0\,.
\end{align}

\vspace{5mm}
We now collect the results established above. We have the following statement.

\vspace{7mm}
{\sc Theorem~4.2:}\,\quad{\em Let the assumptions} {\bf (A1)}--{\bf (A5)} {\em be satisfied, and let 
$\,\bar u\in\uad$ be an optimal control for $({\bf P}_0)$ with the
associated solution $(\bar\mu,\bar\rho,\bar\xi)$ to the state system} (\ref{ss1})--\pier{(\ref{ss5})} {\em in 
the sense of Theorem 2.1. Moreover, let $\{\alpha_n\}\subset (0,1]$ with $\alpha_n\searrow 0$ as
$n\to\infty$ be such that there are optimal pairs $((\bar\mu^{\alpha_n},\bar\rho^{\alpha_n}),\bar u
^{\alpha_n})$ for the adapted problem $(\widetilde{\bf P}_{\alpha_n})$ satisfying} \eqref{3.10}--\eqref{3.12}
{\em (such sequences exist by virtue of Theorem 3.4) and having the associated adjoint variables $\{(p^{\alpha_n},
q^{\alpha_n})\}$. 
Then, for any subsequence 
$\{n_k\}_{k\in\nz}$ of $\nz$, there are a subsequence $\,\{n_{k_\ell}\}_{\ell\in\nz}\,$ and some
triple $(p,q,\lambda)$ such that }
\begin{itemize}
\item $p\in \gianni{H^1(0,T;H)\cap L^\infty(0,T;V)}\cap L^2(0,T;W)$, $q\in L^\infty(0,T;H)$, and $\lambda\in Y'$,
\item {\em the relations} (\ref{eq:4.35})--(\ref{eq:4.38}), (\ref{eq:4.45}), {\em and} (\ref{eq:4.46})
{\em are valid (where the sequences are indexed by $\,n_{k_\ell}\,$ and the limits are taken for
$\ell\to\infty$), and}  
\item {\em the variational inequality} (\ref{vug2}) {\em and the adjoint system equations} \eqref{nadj1},
\eqref{nadj2} and (\ref{eq:4.44})
{\em are satisfied, where \eqref{eq:4.44} holds for every $v\in Y$.}
\end{itemize}

\vspace{3mm}
{\sc Remark:}\quad\,We are unable to show that the limit triple $(p,q,\lambda)$ solving
the adjoint problem associated with the optimal pair $((\bar\mu,\bar\rho),\bar u)$ is uniquely determined. 
Therefore, it may well happen that the limiting pairs differ for different
subsequences. However, it follows from the variational inequality (\ref{vug2}) 
that, for any such limit, it holds with the orthogonal projection  ${\rm I\!P}_{\uad}$ 
onto $\uad$ with respect to the standard inner product in $\pier{L^2(Q)}$ that in the case $\beta_3>0$ we have 
$\,\bar u\,=\,{\rm I\!P}_{\uad}\left(-\beta_3^{-1}p\right)$.

 
\end{document}